\magnification =1100
\input amstex
\documentstyle{amsppt}
\NoBlackBoxes

\def\rto{--\rightarrow}
\def\bbR{{\Bbb R}}

\def\bbC{{\Bbb C}}

\def\bbP{{\Bbb P}}

\def\calL{{\Cal L}}

\def\calN{{\Cal N}}
\def\calO{{\Cal O}}
\def\calP{{\Cal P}}

\def\calA{{\Cal A}}

\def\calJ{{\Cal J}}

\bigskip
\topmatter
\title 
Polar Cremona transformations
\endtitle

\author 
Igor V. Dolgachev 
\endauthor

\thanks 
Research partially supported by NSF Grant DMS 99-70460 and the 
Clay Mathematical Institute
\endthanks

\endtopmatter

\document

Let $F(x_0,\ldots,x_n)$ be a complex homogeneous polynomial of degree $d$.
Consider the linear system
$\calP_F$ generated by the partials
${\partial F\over
\partial x_i}$. We call it the {\it polar linear system} associated to $F$. The problem is to
describe those
$F$ for which the polar linear system is homaloidal, i.e.  the
map 
$(t_0,\ldots,t_n)\to ({\partial F\over
\partial x_0}(t),\ldots,{\partial F\over
\partial x_n}(t))$ is
a birational map. We shall call $F$ with such property a {\it homaloidal
polynomial}. In this paper we review some 
known results about homaloidal polynomials and
also classify them in the cases when $F$ has
no multiple factors and either $n = 3$ or $n = 4$ and
$F$ is the product of linear polynomials.

I am grateful to Pavel Etingof, David Kazhdan and Alexander Polishchuk for
bringing to my attention the problem of classification of homaloidal
polynomials and for various conversations on this matter. Also I thank Hal
Schenck for making useful comments on my paper.

\head {\bf 1. Examples}\endhead

As was probably first noticed by L. Ein and N. Shepherd-Barron 
(see [ES]), many examples of homaloidal polynomials arise from the theory of
prehomogeneous vector spaces. Recall that a complex vector space $V$ is
called  {\it prehomogeneous} with respect to a linear rational representation
of an algebraic group
$G$ in $V$ if there exists a non-constant polynomial $F$ such that the complement of its
set of zeroes is homogeneous with respect to $G$. The polynomial $F$ is necessarily
homogeneous and an eigenvector for $G$ with some character $\chi:G\to \text{GL}(1)$. It generates
the algebra of invariants for the group $G_0 = Ker(\chi)$. The reduced part $F_{red}$ of $F$
(i.e. the product of irreducible factors of $F$) is determined uniquely up to a scalar
multiple.  A prehomogeneous space is called regular if  the determinant of the Hesian
matrix of
$F$ is not identically zero. This definition does not depend on the choice of $F$. We shall
call $F$ a {\it relative invariant} of
$V$.  Note that there is a complete
classification  of regular irreducible prehomogeneous spaces with respect
to a reductive group $G$  (see [KS]).

\medskip
\plainproclaim Theorem ([ES, EKP])]) 1. Let $V$ be a regular
prehomogeneous vector space. Then its relative invariant is a homaloidal
polynomial.

\bigskip
Here are some examples: 
\medskip\noindent
{\bf Examples 1-4.} 1. Any non-degenerate quadratic form $Q$ is obviously a
homaloidal polynomial. The corresponding birational map is a projective
automorphism. It is also a relative invariant for the orthogonal group
$O(Q)\times \text{GL}(1)$ in its natural linear representation. 

2. A reduced cubic polynomial $F$ on $V$ is a relative invariant for a
regular prehomogeneous space with respect to a reductive group $G$ if and only
if the pair
$(V,G)$ is one of the following (up to a linear transformation):

2.1: $G = \text{GL}(1)^3\subset \text{GL}(3), V = \bbC^3$, the action is natural,  $F =
x_0x_1x_2$.

2.2: $G = \text{GL}(3)$, $V $ is the space of quadratic forms on $\bbC^3$, the
action is via the natural action on $\bbC^3$, $F$ is the discriminant
function.

2.3: $G = \text{GL}(3)\times \text{GL}(3)$, $V = \text{Mat}_3$ is the space of
complex $3\times 3$-matrices, the action is by $(g,g')\cdot A = gAg'{}^{-1}$, the
polynomial $F$ is the determinant. 

2.4: $G = \text{GL}(6), V = \Lambda^2(\bbC^6)$, the action is via the natural
action on $\bbC^6$. The polynomial $F$ is the pffafian polynomial.

2.5: $G = E_6\times \text{GL}(1)$, $V = \bbC^{27} = \text{Mat}_3\times \text{Mat}_3\times \text{Mat}_3$ is its
irreducible representation of minimal dimension. The polynomial $F$ is the
Cartan cubic $F(A,B,C) = |A|+|B|+|C|-\text{Tr}(ABC)$. 

The last four examples correspond to the four Severi varieties: 
nonsingular nondegenerate 
subvarieties $S$ of $\bbP^r$ of dimension ${2r-4\over 3}$ whose secant
variety $\text{Sec}(S)$ is not equal to the whole space. The zero locus of the cubic
$F$ in $\bbP(V)$ defines the secant variety. The singular locus of 
$\text{Sec}(S)$
is the Severi variety. According to a theorem from [ES], any homaloidal
cubic polynomial $F$ such that the singular locus of $F^{-1}(0)$ in
$\bbP(V)$ is nonsingular coincides with one from examples 2.2-2.5.

3. Let us identify $\bbP^{n^2-1}$ with the space $\bbP(\text{Mat}_n)$. The map 
$A\to A^{-1}$ is obbviously birational. It is given by the polar linear
system of the polynomial $A\to \det(A)$. The polynomial is a relative
invariant from Example 2.3 (extended to any dimension).

4. The polynomial $F = x_0(x_0x_2+x_1^2)$ is homaloidal. It is a relative
invariant for a prehomogeneous space with respect to a non-reductive group.

\head {\bf 2. Multiplicative Legendre transform}\endhead 

This section is almost entirely borrowed from [EKP]. Let $F\in
\text{Pol}_d(V)$ be a homogeneous polynomial of degree $d$ on a complex vector
space
$V$ of dimension $n+1$.  We denote by $F'$ or by $dF$ the derivative map
$V\to V^*, v\to (dF)_v$. If no confusion arises we also use this notation
for the associated rational map
$\bbP(V)\to
\bbP(V^*)$. If we choose a basis in $V$ and the corresponding dual basis in
$V^*$, we will be able to identify both spaces with $\bbC^n$, and the map
$F'$ with the polar map defined in the introduction. Suppose $F$ is
homaloidal, i.e.
$F'$ defines a birational map
$\bbP(V)\to
\bbP(V^*)$. Then, obviously, 
$d\ln F = F'/F$ defines a birational map $V\to V^*$. 

\plainproclaim Lemma 1. Let $f$ be a homogeneous function of degree $k$ on
$V$ (defined on an open subset) such that $\det (\text{Hess}(\ln f))$ is not
identical zero. Then there exists a homogeneous function
$f_*$ on $V^*$ of degree $k$ such that on some open subset of $V$
$$f_*(d\ln f) = 1/f.\eqno (2.1)$$

{\sl Proof.} Recall first the definition of the {\it Legendre transform}. Let
$Q$ be a  function on $V$ defined in an open neighborhood of a point
$v_0$ such that $\det \text{Hess}(Q)(v_0)\ne 0$. Let $dQ(v_0) = p_0\in V^*$.
Then the Legendre transform $L(Q)$ of $Q$ is the function $L(Q)$ on
$V^*$ defined in an neighborhood of $p_0$ such that 
$$L(Q)(p) = p(v_p)-Q(v_p),\eqno (2.2)$$
where $v_p$ is the unique critical point of the function $v\to p(v)-Q(v)$ in
a neighborhood of $v_0$.

Since the critical point $v_p$ satisfies $p= dQ(v_p)$, we obtain from (2.2)
the equality of functions on an neighborhood of $v_p$ in $V$
$$L(Q)(dQ(v)) = dQ(v)(v) - Q(v).$$
Now let us apply this to $Q = \ln f$. We have
$$L(\ln f)(d\ln f(v)) = d\ln f(v)\cdot v -\ln f(v).$$
Recall that a homogeneous function $H$ of degree $k$ satisfies the Euler
formula:
$$kH(v) = dH(v).$$
Applying this to $H = \ln f$, we get
$$e^{L(\ln f)-k}(d\ln f) = 1/f.$$
It remains to define $f_*$ by 
$$\ln f_* = L(\ln f)-k.\eqno (2.3)$$
It is immediately checked that it is homogeneous of degree $k$.
\bigskip
The function $f_*$ is called the {\it multiplicative Legendre transform of
$f$}. 

\bigskip
\plainproclaim Theorem 2 ([EKP]). Let $F\in \text{Pol}_d(V)$ such that $\det \text{Hess}(\ln
F)$ is not identical zero.  Then $F$ is homaloidal if and only if its
multiplicative transform $F_*$ is a rational function. Moreover, in this case
$$d\ln F_* = (d\ln F)^{-1}.\eqno (2.4)$$

{\sl Proof.} Suppose $F$ is homaloidal. Then $d\ln F$ is a rational
 map of topological degree 1 in its set of definition. It follows from the
definition of the Legendre transform that
$L(\ln F)$ is one-valued on its set of definition. Differentiating (2.1) we
obtain $(d\ln F_*)\circ (d\log F) = \text{id}.$ This checks (2.4). Since 
$d\ln F_* = dF_*/F$ is a homogeneous rational function, the function 
$F_*$ 
must be rational. Conversely, if $F_*$ is rational, we get (2.4) locally, by
differentiating (1). Since $d\ln F_*$ is rational, we have (2.4) globally, and
hence $d\ln F$ is invertible. This implies that $dF$ defines a birational
map, and hence $F$ is homaloidal.

\bigskip
\plainproclaim Corollary 1. Let $F(x_0,\ldots,x_n)$ be a  homaloidal
polynomial of degree
$k > 2$. Assume $F_*$ is a reduced polynomial. Then 
$$k|2(n+1).$$

{\sl Proof.} By the previous theorem 
$$dF_*\circ dF = F^{k-1}(x)F_*(x)(x_0,\ldots,x_n).$$
This implies that the image of the hypersurface $F = 0$ under the birational map $dF:\bbP^n\to
\bbP^n$ is contained in the set of base points of the polar linear system of $F_*$. Since $F_*$
is reduced the latter is a closed subset of codimension $> 1$. Thus $F=0$ is 
contained in the
set of critical points of $dF$ (considered as a map of vector
spaces) and hence $F$ divides the
Hessian determinant. The assertion follows from this.
 
\bigskip
A natural question posed in [EKP] is the following: For which homogenous
polynomials
$F$ its multiplicative Legendre transform
$F_*$ is a polynomial function? 

A  polynomial with this property will be called a {\it homaloidal
EKP-polynomial}.  It is easy to see that $F_*$ has the same degree as
$F$ and $(F_*)_* = F$. It is conjectured  that any homaloidal
EKP-polynomial is a relative invariant of a regular prehomogeneous space.
The converse is proved in [EKP]. In this case $F_* = F$, up to a scaling.

 A remarkable result of [EKP] is the
following:

\plainproclaim Theorem 3. A homaloidal EKP-polynomial of degree 3
coincides with one from Examples 2.

\medskip\noindent
{\bf Example 5.} Consider the polynomial $F$ from Example 4. We have
$$d\ln F = \big({2x_0x_2+x_1^2\over x_0(x_0x_2+x_1^2)},{2x_1\over
x_0x_2+x_1^2},{x_0\over x_0x_2+x_1^2}\big) .$$
Inverting this map we obtain
$$(d\ln F)^{-1} = \big({8x_2\over 4x_0x_2+x_1^2},{4x_1\over 4x_0x_2+x_1^2},
{4x_0x_2-x_1^2\over (4x_0x_2+x_1^2)x_2}\big) = d\ln {(4x_0x_2+x_1^2)^2\over
x_2}.$$ Thus the multiplicative Legendre transform of $F$ equals
$$F_* = {(4x_0x_2+x_1^2)^2\over x_2}.$$
It is a homogeneous rational but not polynomial function.

\head {\bf 3. Plane polar Cremona transformations}\endhead

Here we shall classify all homaloidal polynomials in three variables with no
multiple factors.
  
Since the set of common zeroes of the polars
$\partial_i F$ is equal to the set of non-smooth points of the subscheme $V(F)$, this is
equivalent to requiring that the polars
$\partial_i F$ have no common factors, i.e. the linear system $\calP_F$ has no
fixed part. 

Let $f:\bbP^2\rto \bbP^2$ be a rational map defined by homogeneous polynomials 
$(P_0,P_1,P_2)$ of degree $d$ without common factors. Let $\calJ(f)
\subset k[x_0,x_1,x_2]$ be the ideal generated by the polynomials
$P_0,P_1,P_2$. The corresponding closed subscheme
$B_f = V(\calJ(f))$ of $\bbP^2$ is the base locus subscheme of the linear system spanned
by
$P_0,P_1,P_2$.
 The quotient sheaf
$\calO_{\bbP^2}/\calJ(f)$ is artinian and we denote by $\tilde\mu_x(f)$ the
length of its stalk at a point $x\in V(\calJ(f))$. 

\medskip
\plainproclaim Lemma 2. 
$$\sum_{x\in \bbP^2}\tilde\mu_x(f) = d^2-d_t,$$
where $d_t$ is the degree of the map $f$.

{\sl Proof.} See [Fu], 4.4.

\bigskip
Recall that for any singular point $x$ of $V(F)$ we have the conductor
invariant $\delta_x$ defined as the length of the quotient module 
$\bar\calO_{C,x}/\calO_{C,x}$, where $\bar\calO_{C,x}$ is the normalization
of the local ring $\calO_{C,x}$. Let $r_x$ denote the number of local
branches of $C$ at $x$. We have the following

\plainproclaim Lemma 3. Let $\tilde\mu_x =\tilde\mu_x(f)$, where $f$ is
the map defined by the polar linear system $\calP_F$. For any
$x\in C$,
$$\tilde\mu_x \le  2\delta_x-r_x+1.\eqno (3.1)$$

{\sl Proof.} Without loss of generality we may assume that $x = (1,0,0)$. Let 
$\tilde P(X,Y)$ denote the dehomogenization of a homogeneous polynomial
$P$ with respect to the variable $x_0$. Applying the Euler formula $dF =
x_0F_0+x_1F_1+x_2F_2$, we obtain that 
$$\calJ_x = (\tilde F,{\partial\tilde F\over \partial
X},{\partial \tilde F\over \partial
Y})_x.$$
By Jung-Milnor's formula (see [Mi], Theorem 10.5), the length $\mu_x$ of the
module
$(k[X,Y]/({\partial \tilde F\over \partial
X},{\partial \tilde F\over \partial
Y}))_x$ is equal to $2\delta_x-r_x+1$. It remains to observe that
$\tilde\mu_x \le \mu_x$.

\medskip 
The next lemma is a well-known formula for the arithmetic genus of a
plane curve.

\plainproclaim Lemma 4. 
$$p_a(C) = (d-1)(d-2)/2 = \sum_{i=1}^hg_i+\sum_x\delta_x-h+1,\eqno (3.2)$$
where $h$ is the number of irreducible components $C_i$ of $C$ and $g_i$ is
the genus of the normalization of $C_i$.

The next formula is an easy consequence of the incidence relation 
count for pairs of lines, 
but just for fun we give a high-brow proof of this:

\medskip
\plainproclaim Corollary 2. Let $\{L_1,\ldots,L_s\}$ be a set of lines in
$\bbP^2$.  Let $a_i$ denote the number of points which belong to $i\ge 2$
distinct lines. Then
$$s(s-1) = \sum_{i=2}^s a_ii(i-1).\eqno (3.3)$$

{\sl Proof.} We apply the previous formula to the curve $L = L_1+\ldots+L_s$.
Each singular point of $L$ lies on the intersection of $i\ge 2$ lines.
It is isomorphic locally to the singular point of the affine curve given by
an equation $\prod_{j=1}^i(\alpha_jX+\beta_jY) = 0$. It is easy to compute
$\delta_x$. It is equal to $i(i-1)/2$. Since $r_x = i$, by Lemma 4, we
have
$$(s-1)(s-2)/2 = \sum_{i=2}^s a_ii(i-1)/2-s+1.$$
This is equivalent to the claimed formula.

\bigskip
\plainproclaim Theorem 4. Let $F$ be a homaloidal polynomial in three variables without 
multiple factors. Then, after a linear change of variables, it coincides with one from Examples
1, 2.1, 4. In other words, 
$C = V(F)$ is one of the following curves:
\item{(i)} a nonsingular conic;
\item {(ii)} the union of three nonconcurrent lines;
\item{(iii)} the union of a conic and its tangent.

{\sl Proof}.  Since $\calP_F$ is
homaloidal, we can apply Lemma 2 to obtain
$$d^2-2d = \sum_{x\in C}\tilde\mu_x.\eqno (3.4)$$
By Lemma 3,
$$d^2-2d \le \sum_{x\in C}(2\delta_x-r_x+1).$$
By Lemma 4,
$$d^2-3d = 2\sum_{i=1}^hg_i+2\sum_{x\in C}\delta_x-2h.\eqno (3.5)$$
Let $C_1,\ldots,C_h$ be irreducible components of $C$ and $d_i = \deg
C_i$. Using (3.4) and (3.5), we obtain
$$\sum_{i=1}^h(2-d_i) = -d+2h \ge 2\sum_{i=1}^hg_i+\sum_{x\in C}(r_x-1)\ge
0.\eqno (3.6)$$

The rest of the proof consists of analyzing this inequality. First observe
that each point of intersection of two irreducible components gives a positive
contribution to the sum 
$\sum_{i=1}^k(r_i-1)$. This immediately implies that $d_i = 1$ for some $i$
unless $C$ is an irreducible conic. In the latter case it is obviously
nonsingular (otherwise the polar linear system is a pencil). This is case
(i) of the theorem. So we may assume that $C_1,\ldots,C_s$ are lines. It
follows from (3.6) that 
$$0\ge \sum_{i=s+1}^h(2-d_i) \ge 2\sum_{i=1}^hg_i+\sum_{x\in C}(r_x-1)-s.\eqno
(3.7)$$ 
If $s = 1$, then each point of intersection of $C_1$ with other
component of $C$ contributes at least 1 to the sum $\sum_{i=1}^k(r_i-1)$.
This shows that $C = C_1+C_2$, where $L$ intersects $C_2$ at one point
and
$d_2 = 2$. This is case (iii) of the theorem. 

Assume that $s\ge 2$. Let
$x_1,\ldots,x_N$ be the intersection points  of the lines
$C_1,\ldots,C_s.$ Let $a_j$ be the number of points among them which
belong to $j\ge 2$ lines. Then
$\sum_{j=2}^sa_j = N$, and
$$\sum_{x\in C}(r_x-1)-s\ge \sum_{i=1}^N(r_i-1)-s\ge
\sum_{j=2}^sja_j-N-s = \sum_{j=2}^s(j-1)a_j-s.\eqno (3.8)$$ By (3.3), 
$$s = \sum_{j=2}^s{j\over  s-1}a_j(j-1).$$
Assume not all lines pass through one point, i.e.  $a_s = 0$. Then
$j\le s-1$ for all $j$ with $a_j\ne 0$. In this case
$$s\le \sum_{j=2}^sa_j(j-1)\eqno (3.9)$$ and the equality holds if and only
if 
$a_j = 0$ for all $j\ne s-1$. If $p_i$ is a point lying on $s-1$ lines, then
the remaining line must intersect other lines at points different from
$p_i$. This gives that $a_2\ne 0$. So, if the equality holds, we have $s = 3$ and $a_2 = N = 3$.
If $h\ne s$,  then $C_h$ is of degree $> 1$. Its points of intersection with
three lines give positive contribution to the sum 
$\sum_{x\ne x_1,\ldots,x_N}(r_x-1)-s$.
Thus (3.8) is a strict inequality contradicting (3.7). So $C$ is the
union of three nonconcurrent lines which is case (ii) of the theorem.

  It remains to consider the case when all lines pass
through one point. In this case  (3.7) implies that $s < h$. Then
 $C_h$ is of degree $ > 1$. Assume $p_1\in C_h$. Then $r_1 \ge s+1$ and
$$\sum_{x\in C}(r_x-1)-s = (r_1-1-s)+\sum_{x\ne p_1}(r_x-1) \ge 0.\eqno (3.9) $$
It follows from (3.7) that $C_h$ is a nonsingular conic. Since
$s\ge 2$, one of the lines is not tangent to $C_h$ at $p_1$ and hence
intersects $C_h$ at some point $x\ne p_1$. Thus (3.9) is a strict inequality.
  This
contradicts (3.7).  If $p_1\not\in C_h$, then $C_h$ intersects each line so
that we have $\sum_{x\ne p_1}(r_x-1)\ge s$ and 
$$\sum_{x\in C}(r_x-1)-s = (r_1-1-s)+\sum_{x\ne p_1}(r_x-1) \ge s-1 > 0. $$ 
Again a contradiction. 
 
\bigskip 
Let us notice the following combinatorial fact which follows from
the proof of the above theorem in the case when $C$ is the union of lines:

\plainproclaim Corollary 3. Let $C$ consist of $s$ lines $l_1,\ldots,l_s$.
For each line $l_i$ let
$k_i$ be the number of singular points of $C$ on $l_i$, and let $t$ be the total number of
singular points. Assume that $t > 1$. Then
$$\sum_{i=1}^s(k_i-1)\ge t$$
and the equality takes place if and only if $t = 3, s = 3$.

{\sl Proof.} Let $d$ be the degree of the map given by the polar linear system of the 
polynomial defining $C$. We resolve the indeterminacy points by blowing up the singular
points of $C$. Let $E_p$ be the exceptional curve blow-up from the point $p$, $h$ be the
class of a general line and $m_p$ be the multiplicity of a singular point $p$. Then 
$$
d = ((s-1)h-\sum_{p\in \text{Sing}(C)}(m_p-1)E_p)^2= (s-1)^2-\sum_{p\in
\text{Sing}(C)}(m_p-1)^2.$$ Let $a_i = \#\{p:m_p =i\}$. Applying
equality (3.3), we can rewrite it as follows:
$$d = s(s-1)-(s-1)-\sum_{i=2}^sa_i(i-1)i+\sum_{i=2}^sa_i(i-1) =$$
$$-(s-1)+\sum_{i=2}^sa_i(i-1) = -s+1+\sum_{i=2}^sia_i-\sum_{i=2}^sa_i.$$
Now the standard incidence relation argument gives us 
$$\sum_{i=2}^sia_i = \sum_{p\in \text{Sing}(C)}m_p = \sum_{i=1}^sk_i.$$
This allows us to rewrite the expression for $d$ in the form
$$d = 1+\sum_{i=1}^s(k_i-1)-t.$$
Now $d\ge 1$ unless all lines pass though one point and, by  Theorem 4,
$d = 1$ if and  only if $s=3,t=3$.

\medskip\noindent
{\bf Remark} As was explained to me by Hal Schenck, in the case of
a real arrangement of lines the previous Corollary follows easily from the
Euler formula applied to the cellular subdivision of $\bbR\bbP^2$ defined
by the arrangement. One interprets the left-hand side as the number $f_1$ 
of edges, the right-hand side as the number $f_0$ of vertices and uses
that 
$f_0\ge s$ and $f_2\ge f_0+1$ if the arrangement is not a pencil (see [Gr],
pp.10 and 12).

 \bigskip
Unfortunately, the argument used in the proof of Theorem 4 does not apply to non-reduced
polynomials. However, the following conjecture seems to be reasonable:

\bigskip
\plainproclaim Conjecture. Let $F = A_1^{m_1}\ldots A_s^{m_s}$ be the
factorization of $F$ into prime factors. Let $G = A_1\ldots
A_s$. Then the polar linear system
$\calP_F$ is homaloidal if and only if
$\calP_{G}$ is homaloidal.

\head {\bf 4. Arrangements of hyperplanes in $\bbP^3$}\endhead

Here we shall consider the special case when $F = \prod_{i=1}^nL_i$ is the product of linear
polynomials in four variables without multiple factors. Its set of zeroes is an arrangement
of hyperplanes in
$\bbP^3$.

Let $\calA = \{H_1,\ldots,H_N\}$ be the set of planes $\{L_i = 0\}$,
$\calL$ be the set of lines which are contained in more than one plane
$H_i$, and $\calP$ be the set of points which are contained in more than
two planes $H_i$. For any $l\in \calL$, set
$$k_l = \#\{i:l\subset H_i\}, \quad a_l = \#\{p\in \calP:p\in l\}.$$
For any $p\in \calP$ set
$$k_p = \#\{i:p\in H_i\}.$$
We define $d_\calA$ to be the degree of the polar linear system defined by $F$.

\plainproclaim Lemma 5.
$$d_\calA  = (N-1)^3-\sum_{p\in \calP}(k_p-1)+\sum_{l\in\calL}(k_l-1)(a_l-1).$$

{\sl Proof.} We can resolve the points of indeterminacy of $\calP_F$ by first blowing up
each points $p\in \calP$ followed by blowing up the proper transforms of each line $l\in
\calL$. Let 
$$D = \sum_{p\in \calP}(k_p-1)E_p+\sum_{l\in\calL}(k_l-1)E_l.$$
Here the notations are self-explanatary. We have (see [Fu])
$$d_\calA = ((N-1)H-D)^3,$$
where $H$ is the pre-image of a general plane in the blow-up. Using the standard formulae
for the blow-up a smooth subvariety, we have
$$E_l^3 = -c_1(N_{\bar l}) = -[(4H-2\sum_{l\in \calL,p\in l}E_p)\cdot\bar l-2] = 2a_l-2.$$
Here $\bar l$ denotes the proper transform of the line $l$ under the blowing up the points
from $\calP$, and $N_{\bar l}$ is the normal bundle of $\bar l$. Next, we have
$$E_l^2\cdot E_p = -1,\quad E_p^3 = 1.$$
Collecting this together we get
$$D^3 = \sum_{l\in\calL}(k_l-1)^3(2a_l-2)+\sum_{p\in\calP}(k_p-1)^3-
3\sum_{l\in \calL,p\in
l}(k_l-1)^2(k_p-1),$$
$$H\cdot D^2 = \sum_{l\in\calL}(k_l-1)^2E_l\cdot H = -\sum_{l\in\calL}(k_l-1)^2,$$
$$H^2\cdot D = 0.$$
This gives
$$d_\calA = (N-1)^3-3(N-1)\sum_{l\in\calL}(k_l-1)^2-\sum_{l\in\calL}(k_l-1)^3(2a_l-2)-$$
$$\sum_{p\in\calP}(k_p-1)^3+3\sum_{l\in \calL,p\in
l}(k_l-1)^2(k_p-1).$$
Observe now that
$$\sum_{p\in l}(k_p-1) = \sum_{p\in l}k_p-a_l = (a_lk_l+N-k_l)-a_l =
(a_l-1)k_l+N-a_l.$$
This allows us to rewrite the expression for $d$ as follows:
$$d_\calA = (N-1)^3-3(N-1)\sum_{l\in\calL}(k_l-1)^2-\sum_{l\in\calL}(k_l-1)^3(2a_l-2)-$$
$$\sum_{p\in\calP}(k_p-1)^3+3\sum_{l\in \calL}(k_l-1)^3(a_l-1)+3(N-1)\sum_{l\in
\calL}(k_l-1)^2 =$$
$$(N-1)^3-\sum_{p\in \calP}(k_p-1)+\sum_{l\in\calL}(k_l-1)(a_l-1).$$
This proves the lemma.

\plainproclaim Lemma 6. Let 
$$t_s = \#\{p\in \calP:k_p = s\}, \quad t_q(1) = \#\{l\in \calL:k_l = q\},$$
$$t_{sq} = \sum_{l\in\calL:k_l=q}\#\{p\in l:k_p = s\}.$$
Then
$$\binom N3 = \sum_{s}\binom s3t_s-\sum_{s,q}\binom q3(t_{sq}-t_q(1)).$$

{\sl Proof.} This is a three-dimensional analog of Corollary 2 to Lemma
4.  It  easily follows from the incidence relation count for triples of
distinct planes and points and lines. 

\plainproclaim Corollary 4. 
$$d_\calA = N-1-\sum_{p\in\calP}(k_p-1)+\sum_{l\in\calL}(a_l-1)(k_l-1).$$

{\sl Proof.} Combine the previous two lemmas.

\plainproclaim Lemma 7. Let $\calA$ be an arrangement of $N$ hyperplanes
in $\bbP^3$ defined by a polynomial $F$. The following properties are
equivalent:
\item{(i)}  all planes pass through a point;
\item{(ii)} the partials of $F$ are linearly dependent
\item{(ii)} $d_\calA = 0$.

{\sl Proof.} Obvious.

\medskip 
\plainproclaim Lemma 8. Let $\calA$ be an arrangement of $N$ planes. Let
$\calA'$ be a  new arrangement obtained by adding one more plane to
$\calA$. Assume $d_\calA \ne 0$. Then
$$d_{\calA'} > d_{\calA}.$$

{\sl Proof.} Let 
$$\calP' = \{p\in\calP:p\in H\}, \quad \calL' =\{l\in \calL:l\subset H\},$$ 
$$ \calL'' = \{l\in \calL:p\not\in l\quad \hbox{for any $p\in \calP'$}\},$$
$$ \calN =
\{l\subset H\cap (H_1\cup\ldots\cup H_N)\}\setminus \calL.$$ 
Note that each line $l\in\calN$ is a double line and each line $l\in \calL''$ contains one
new singular point $H\cap l$ of multiplicity $k_l+1$. Applying the 
previous corollary, we
obtain
$$d_{\calA'} =  N-\sum_{p\in \calP\setminus\calP'}(k_p-1)-\sum_{p\in \calP'}k_p-\sum_{l\in\calL''}k_l+\sum_{l\in
\calL'}k_l(a_l-1)+$$
$$\sum_{l\in
\calL\setminus\calL'}(k_l-1)a_l+\sum_{l\in\calN}(a_l'-1),$$
where $a_{l'}$ denotes the number $a_l$ defined for the extended arrangement.
Applying the corollary again, we get
$$d_{\calA'}-d_{\calA} =
1+(\sum_{l\in\calL\setminus(\calL'\cup\calL''}(k_l-1)-\#\calP')+
(\sum_{l\in\calN}(a_l'-1)-\#\calL'')+\sum_{l\in\calL'}(a_l-1).\eqno (4.1)$$
For each $p\in \calP'$ there exists a line $l\in \calL\setminus(\calL'\cup\calL'')$ passing through $p$. Since $k_l > 1$ for each line we see that 
$\sum_{l\in\calL\setminus(\calL'\cup\calL''}(k_l-1)-\#\calP'
\ge 0$. Now consider the arrangement of lines in the plane $H$ formed by the lines
$l\in\calN$. Its multiple points are the points of intersection of $H$ with lines in
$\calL''$. Applying Corollary 3 to Theorem 4, we see that 
$\sum_{l\in\calN}(a_l'-1)-\#\calL''\ge 0$ unless there is only one line in $\calL''$ when
this difference is equal to $-1$. But in this case $H$ must contain at least one line from
$\calL$ and hence there is an additional term $\sum_{l\in\calL'}(a_l-1)$. If it is zero,
then each line $l\in\calL'$ contains only one singular point of the arrangement. This
implies that all planes except maybe one contain $l$. In this case all planes pass through
a point and $d_\calA = 0$. So the term is positive and we have proved the inequality
$d_{\calA'} > d_{\calA}$.

\plainproclaim Theorem 5. Let $\calA$ be an arrangement of $N$ planes in $\bbP^3$ with
$d_\calA = 1$. Then $\calA$ is the union of four planes in general linear position.

{\sl Proof.} By the previous lemma deleting any plane $H$ from the arrangement $\calA$ defines
an arrangement $\calA'$ with $d_{\calA} = 0$. 
We may assume that $H$ does not pass through the common point of the planes from
$\calA'$. In the notation of
the proof of the previous lemma, where the new arrangement
 is our $\calA$ and the
old one is $\calA\setminus \{H\}$, we have $\#\calL'' = N-1$. Now the term 
$(\sum_{l\in\calN}(a_l'-1)-\#\calL'')$ in (4.1) must be equal to zero since otherwise
$d_{\calA} > 1$. By Lemma 6, $N-1 = 3$. Thus $N = 4$. Since $d_\calA \ne
0$, the planes do not have a common point and hence the arrangement is as
in the assertion of the theorem.

\bye